\documentstyle[twoside,11pt,leqno]{article}

\def\C{{\bf C}}
\def\ind{\textrm{ind}}
\def\X{{\mathcal X}}
\def\Y{{\mathcal Y}}

\def\C{\texttt{C}}
\def\iso{\textrm{iso}}

\def\a{{\Xi(A)}}
\def\b{{\Xi(B)}}
\def\A{{\Xi(A^*)}}
\def\B{{\Xi(B^*)}}
\def\M{{M_C}}
\def\m{{M_0}}
\def\asc{ \textrm{asc}}
\def\dsc{ \textrm{dsc}}

\newtheorem{df}{Definition}[section]
\newtheorem{thm}[df]{Theorem}
\newtheorem{prop}[df]{Proposition}
\newtheorem{cor}[df]{Corollary}
\newtheorem{ex}[df] {Example}
\newtheorem{rem}[df] {Remark}

\def\subject#1{\renewcommand{\thefootnote}{}
\footnote{AMS(MOS) subject classification (2000). {#1}}}

\def\keywords#1{\renewcommand{\thefootnote}{}
\footnote{ Key words and phrases. {#1}}}

\expandafter\let\csname enddemo*\endcsname=\enddemo

\def\qedsymbol{\ifmmode\bgroup\else$\bgroup\aftergroup$\fi
  \vcenter{\hrule\hbox{\vrule
height.5em\kern.5em\vrule}\hrule}\egroup}
\def\qed{\ifmmode\else\unskip\nobreak\fi\quad\qedsymbol}
\def\wt#1{\widetilde{#1}}

\title
{\bf  Upper Triangular Operator Matrices, SVEP and Browder, Weyl
Theorems \/}

\author {\normalsize B. P. Duggal } \vskip 1truecm

\date{}

\begin{document}

\maketitle \thispagestyle{empty} %\baselineskip=14pt

\subject{Primary 47B47, 47A10, 47A11} \keywords{  Banach space,
 operator matrix,  Browder spectrum,  Weyl spectrum, single valued extension
 property,polaroid operator
 Browder and  Weyl theorems.}

\vskip 1truecm

\setlength{\baselineskip}{14pt}

\date{}

\begin{abstract} A Banach space operator $T\in B(\X)$ is polaroid if
points $\lambda\in\iso\sigma\sigma(T)$ are poles of the resolvent
of $T$. Let $\sigma_a(T)$, $\sigma_w(T)$, $\sigma_{aw}(T)$,
$\sigma_{SF_+}(T)$ and $\sigma_{SF_-}(T)$ denote, respectively,
the approximate point, the Weyl, the Weyl essential approximate,
the upper semi--Fredholm and lower semi--Fredholm spectrum of
$T$. For $A$, $B$ and $C\in B(\X)$, let $M_C$ denote the operator
matrix $\left(\begin{array}{clcr} A & C\\0 &
B\end{array}\right)$. If $A$ is polaroid on
$\pi_0(\M)=\{\lambda\in\iso\sigma(\M):
0<\dim(\M-\lambda)^{-1}(0)<\infty\}$, $\m$ satisfies Weyl's
theorem, and $A$ and $B$ satisfy either of the hypotheses (i) $A$
has SVEP at points
$\lambda\in\sigma_w(\m)\setminus\sigma_{SF_+}(A)$ and $B$ has
SVEP at points $\mu\in\sigma_w(\m)\setminus\sigma_{SF_-}(B)$, or,
(ii) both $A$ and $A^*$ have SVEP at points
$\lambda\in\sigma_w(\m)\setminus\sigma_{SF_+}(A)$, or, (iii) $A^*$
has SVEP at points
$\lambda\in\sigma_w(\m)\setminus\sigma_{SF_+}(A)$ and $B^*$ has
SVEP at points $\mu\in\sigma_w(\m)\setminus\sigma_{SF_-}(B)$, then
$\sigma(\M)\setminus\sigma_w(\M)=\pi_0(\M)$. Here the hypothesis
that  $\lambda\in\pi_0(\M)$ are poles of the resolvent of $A$ can
not be replaced by the hypothesis $\lambda\in\pi_0(A)$ are poles
of the resolvent of $A$.

For an operator $T\in B(\X)$, let
$\pi_0^a(T)=\{\lambda:\lambda\in\iso\sigma_a(T),
0<\dim(T-\lambda)^{-1}(0)<\infty\}$. We prove that if $A^*$ and
$B^*$ have SVEP, $A$ is polaroid on $\pi_0^a(\M)$ and $B$ is
polaroid on $\pi_0^a(B)$ , then
$\sigma_a(\M)\setminus\sigma_{aw}(\M)=\pi_0^a(\M)$.
\end{abstract}

%\numberwithin{equation}{section}

\def\demo{
  \par\topsep5pt plus5pt
  \trivlist
  \item[\hskip\labelsep\it Proof.]\ignorespaces}
\def\enddemo{\qed \endtrivlist}
\expandafter\let\csname enddemo*\endcsname=\enddemo

\def\qedsymbol{\ifmmode\bgroup\else$\bgroup\aftergroup$\fi
  \vcenter{\hrule\hbox{\vrule
height.5em\kern.5em\vrule}\hrule}\egroup}
\def\qed{\ifmmode\else\unskip\nobreak\fi\quad\qedsymbol}
\def\wt#1{\widetilde{#1}}

%%%%%%%%%%%%%%%%%%%%%%%%%%%%%%%%%%%%%%%%%%%%%%%%%%%%%%%%%%%%%%%%%%%%%%%%%%%
%%%%%%%%%%%%%%%%%%%%%%%%%%%%%%%%%%%%%%%%%%%%%%%%%%%%%%%%%%%%%%%%%%%%%%%%%%%%

\section { Introduction} A Banach space
operator $A$, $A\in B(\X)$, is upper semi-Fredholm (resp., lower
semi-Fredholm) at a complex number $\lambda\in\C$ if the range
$(A-\lambda)\X$ is closed and
$\alpha(A-\lambda)=\dim(A-\lambda)^{-1}(0)<\infty$ (resp.,
$\beta(A-\lambda)=\dim(\X/(A-\lambda)\X)<\infty$). Let
$\lambda\in\Phi_+(A)$ (resp., $\lambda\in\Phi_-(A)$) denote that
$A$ is upper semi-Fredholm (resp., lower semi-Fredholm) at
$\lambda$. The operator $A$ is Fredholm at $\lambda$, denoted
$\lambda\in\Phi(A)$, if $\lambda\in\Phi_+(A)\cap\Phi_-(A)$. $A$
is Browder (resp., Weyl) at $\lambda$ if $\lambda\in\Phi(A)$ and
$\asc(A-\lambda)=\dsc(A-\lambda)<\infty$ (resp., if
$\lambda\in\Phi(A)$ and $\ind(A-\lambda)=0$). Here
$\ind(A-\lambda)=\alpha(A-\lambda)-\beta(A-\lambda)$ denotes the
Fredholm index of $A-\lambda$, $\asc(A-\lambda)$ denotes the
ascent of $A-\lambda$ ($=$ the least non-negative integer $n$
such that $(A-\lambda)^{-n}(0)=(A-\lambda)^{-(n+1)}(0)$) and
$\dsc(A-\lambda)$ denotes the descent of $A-\lambda$ ($=$ the
least non-negative integer $n$ such that
$(A-\lambda)^n\X=(A-\lambda)^{n+1}\X$). Let $\sigma(A)$ denote
the spectrum, $\sigma_a(A)$ the approximate point spectrum,
$\iso\sigma(A)$ the set of isolated points of $\sigma(A)$,
$\pi_0(A)=\{\lambda\in\iso\sigma(A):
0<\alpha(A-\lambda)<\infty\}$,
$\pi_0^a(A)=\{\lambda\in\iso\sigma_a(A):
0<\alpha(A-\lambda)<\infty\}$, and $p_0(A)$ the set of finite rank
poles (of the resolvent) of $A$. The Browder spectrum
$\sigma_b(A)$ of $A$ is the set $\{\lambda\in\C: A-\lambda$ is
not Browder$\}$, the Weyl spectrum $\sigma_w(A)$ of $A$ is the
set $\{\lambda\in\C: A-\lambda$ is not Weyl$\}$, the Browder
essential approximate spectrum $\sigma_{ab}(A)$ of $A$ is the set
$\{\lambda\in\C: \lambda\notin\Phi_+(A)$ or $
\asc(A-\lambda)\not<\infty\}$, and the Weyl essential approximate
spectrum $\sigma_{aw}(A)$ of $A$ is the set $\{\lambda\in\C:
\lambda\notin\Phi_+(A)$ or $ \ind(A-\lambda)\not\leq 0\}$.
Following current terminology, the operator $A$ satisfies:
Browder's theorem, or $Bt$, if $\sigma_w(A)=\sigma_b(A)$
(equivalently, $\sigma(A)\setminus\sigma_w(A)=p_0(A)$); Weyl's
theorem, or $Wt$, if $\sigma(A)\setminus\sigma_w(A)=\pi_0(A)$;
$a$--Browder's theorem, or $a-Bt$, if
$\sigma_{aw}(A)=\sigma_{ab}(A)$; $a$--Weyl's theorem, or $a-Wt$,
if $\sigma_a(A)\setminus\sigma_{aw}(A)=\pi_0^a(A)$.

\

 An  operator  $A\in B(\X)$  has {\it the single-valued extension property} at
$\lambda_0\in {\C}$, SVEP at $\lambda_0$, if for every open disc
${\bf {\mathcal D}}_{\lambda_0}$ centered at $\lambda_0$ the only
analytic function $f:{\bf{\mathcal D}}_{\lambda_0}\rightarrow
{\X}$ which satisfies
$$
(A-\lambda)f(\lambda)=0 \hspace{3mm}\mbox{for all} \hspace{3mm}
\lambda\in {\bf{\mathcal D}}_{\lambda_0}
$$
is the function $f\equiv 0$. Trivially, every operator $A$ has
SVEP on the resolvent set $\rho(A)={\C}\setminus \sigma(A)$; also
$A$ has  SVEP at points $\lambda\in \textrm{iso}\sigma(A)$. Let
$\a$ denote the set of $\lambda\in\C$ where $A$ does not have
SVEP: {\em we say that $A$ has SVEP if  $\a=\emptyset$.} SVEP
plays an important role in determining the relationship between
the  Browder and Weyl spectra, and the Browder and Weyl theorems.
Thus $\sigma_b(A)=\sigma_w(A)\cup\a=\sigma_w(A)\cup\A$, and if
$A^*$ has SVEP then $\sigma_b(A)=\sigma_w(A)=\sigma_{ab}(A)$
\cite[pp 141-142]{A}; $A$ satisfies $Bt$ (resp., $a-Bt$) if and
only if $A$ has SVEP at $\lambda\notin\sigma_w(A)$ (resp.,
$\lambda\notin\sigma_{aw}(A)$) \cite[Lemma 2.18]{Du}; and if
$A^*$ has SVEP, then $A$ satisfies $Wt$ implies $A$ satisfies
$a-Wt$ \cite[Theorem 3.108]{A}.

\

For $A$, $B$ and $C\in B(\X)$, let $\M$ denote the upper
triangular operator matrix $\M=\left(\begin{array}{clcr} A & C\\0
& B\end{array}\right)$. A study of the spectrum, the Browder and
Weyl spectra, and the Browder and Weyl theorems for the operator
$\M$, and the related diagonal operator $\m=A\oplus B$, has been
carried by a number of authors in the recent past (see
\cite{{Cao},{DDj},{DZ},{Lee}} for further references). Thus, if
either $\A=\emptyset$ or $\b=\emptyset$, then
$\sigma(\M)=\sigma(\m)=\sigma(A)\cup\sigma(B)$; if
$\a\cup\b=\emptyset$, then $\M$ has SVEP,
$\sigma_b(\M)=\sigma_w(\M)=\sigma_b(\m)=\sigma_w(\m)$, and $\M$
satisfies $a-Bt$. Browder's theorem, much less Weyl's theorem,
does not transfer from individual operators to direct sums: for
example, the forward unilateral shift and the backward unilateral
shift on a Hilbert space satisfy $Bt$, but their direct sum does
not. However, if $(\a\cap\B)\cup\A=\emptyset$, then : $\m$
satisfies $Bt$ (resp., $a-Bt$) implies $\M$ satisfies $Bt$ (resp.,
$a-Bt$); if  points $\lambda\in\iso\sigma(A)$ are eigenvalues of
$A$, $A$ satisfies $Wt$, then $\m$ satisfies $Wt$ implies $\M$
satisfies $Wt$ \cite[Proposition 4.1 and Theorem 4.2]{DZ}. Our
aim in this paper is to fine tune some of the extant results to
prove that: $\sigma_b(\m)=\sigma_b(\M)\cup\{\A\cap\b\}$;
$\sigma_{ab}(\M)\subseteq\sigma_{ab}(\m)\subseteq\sigma_{ab}(\M)\cup\Xi^*_+(A)\cup\Xi_+(B)$;
and
$\sigma_w(A)\cup\sigma_w(B)\subseteq\sigma_w(\M)\cup\{\Xi)P)\cup\Xi(Q)\}$,
where,except for $P=A$ and $Q=B^*$, $P=A$ or $A^*$ and $Q=B$ or
$B^*$, $\Xi_+^*(A)=\{\lambda:\lambda\notin\Phi_+(A)$ or $A^*$
does not have SVEP at $\lambda\}$ and
$\Xi_+(B)=\{\lambda:\lambda\notin\Phi_+(B)$ or $B$ does not have
SVEP at $\lambda\}$. Let $\sigma_{SF_+}(A)$ (resp.,
$\sigma_{SF_-}(A)$) denote the {\em upper semi--Fredholm
spectrum} (resp., {\em lower semi--Fredholm spectrum}) of $A$. It
is proved that if points $\lambda\in\pi_0(\M)$ are poles of $A$,
$\m$ satisfies $Wt$, and $A$ and $B$ satisfy either of the
hypotheses (i) $A$ has SVEP at points
$\lambda\in\sigma_w(\m)\setminus\sigma_{SF_+}(A)$ and $B$ has
SVEP at points $\mu\in\sigma_w(\m)\setminus\sigma_{SF_-}(B)$ or
both $A$ and $A^*$ have SVEP at points
$\lambda\in\sigma_w(\m)\setminus\sigma_{SF_+}(A)$ or $A^*$ has
SVEP at points $\lambda\in\sigma_w(\m)\setminus\sigma_{SF_+}(A)$
and $B^*$ has SVEP at points
$\mu\in\sigma_w(\m)\setminus\sigma_{SF_-}(B)$, then $\M$
satisfies $Wt$. Here the hypothesis that points
$\lambda\in\pi_0(\M)$ are poles of $A$ is essential. We prove also
that if $\A\cup\B=\emptyset$, points $\lambda\in\pi_0^a(\M)$ are
poles of $A$ and points $\mu\in\pi_0^a(B)$ are poles of $B$, then
$\M$ satisfies $a-Wt$.

\

Throughout the following, the operators $A$, $B$ and $C$ shall be
as in the operator matrix $\M$; we shall write $T\in B(\Y)$ for a
general Banach space operator.

\section{ Browder, Weyl Spectra and SVEP}
We start by recalling some results which will be used in the
sequel without further reference.

For an operator $T\in B(\Y)$ such that $\lambda\in\Phi_{\pm}(T)$,
the following statements are equivalent \cite[Theorems 3.16 and
3.17]{A}:

(a) $T$ (resp., $T^*$) has SVEP at $\lambda$;

(b) $\asc(T-\lambda)<\infty$ (resp., $\dsc(T-\lambda)<\infty$).

Furthermore, if $\lambda\in\Phi_{\pm}(T)$, and both $T$ and $T^*$
have SVEP at $\lambda$, then
$\asc(T-\lambda)=\dsc(T-\lambda)<\infty$,
$\lambda\in\iso\sigma(T)$ and $\lambda$ is a pole of (the
resolvent of) $T$ \cite[Corollary 3.21]{A}. For an operator $T\in
B(\Y)$ such that $\lambda\in\Phi(T)$ with $\ind(T-\lambda)=0$,
$T$ or $T^*$ has SVEP at $\lambda$ if and only if
$\asc(T-\lambda)=\dsc(T-\lambda)<\infty$. Evidently,
$\asc(\M)<\infty\Longrightarrow \asc(A)<\infty$. It is not
difficult to verify, \cite[Lemma 2.1]{DZ}, that
$\dsc(B)=\infty\Longrightarrow \dsc(\M)=\infty$. In general,
\begin{eqnarray*} &  &
\sigma(\M)\subseteq\sigma(A)\cup\sigma(B)=\sigma(\m)=\sigma(\M)\cup\{\A\cup\b\};\\
&  & \sigma_b(\M)\subseteq
\sigma_b(A)\cup\sigma_b(B)=\sigma_b(\m);\end{eqnarray*} and
\begin{eqnarray*}
\sigma_w(\M)\subseteq\sigma_w(\m)\subseteq\sigma_w(A)\cup\sigma_w(B).
\end{eqnarray*} If $\lambda\notin\sigma_{aw}(\M)$, then
$\lambda\in\Phi_+(A)$, and either $\alpha(B-\lambda)<\infty$ and
$\ind(A-\lambda)+\ind(B-\lambda)\leq 0$, or
$\beta(A-\lambda)=\alpha(B-\lambda)=\infty$ and $(B-\lambda){\X}$
is closed, or $\beta(A-\lambda)=\infty$ and $(B-\lambda){\X}$ is
not closed \cite[Theorem 4.6]{DDj}. Since $$\M
-\lambda=\left(\begin{array}{clcr} I & 0\\0 &
B-\lambda\end{array}\right)\left(\begin{array}{clcr} I & C\\0 &
I\end{array}\right)\left(\begin{array}{clcr} A-\lambda & 0\\0 &
I\end{array}\right),$$ $\lambda\in\Phi(\M)$ implies that
$\lambda\in\Phi_+(A)\cap\Phi_-(B)$. Using this it is seen that if
$\lambda\notin\sigma_b(\M)$, then
$\lambda\in\Phi_+(A)\cap\Phi_-(B)$,
$\ind(A-\lambda)+\ind(B-\lambda)=0$, $\asc(A-\lambda)<\infty$ and
$\dsc(B-\lambda)<\infty$. If in addition
$\lambda\notin\Xi(P)\cup\Xi(Q)$, where (except for $P=A$ and
$Q=B^*$) $P=A$ or $A^*$ and $Q=B$ or $B^*$, then it is seen (argue
as in the proof of Proposition 2.1 below) that
$\ind(A-\lambda)=\ind(B-\lambda)=0$. Thus
$\lambda\notin\sigma_b(\m)$, which implies that
$$\sigma_b(\m)\subseteq\sigma_b(\M)\cup\{\Xi(P)\cup\Xi(Q)\},$$
where, except for $P=A$ and $Q=B^*$, $P=A$ or $A^*$ and $Q=B$ or
$B^*$. The following proposition gives
more.\begin{prop}\label{prop1}
$\sigma_b(\m)=\sigma_b(\M)\cup\{\A\cap\b\}$.\end{prop}\begin{proof}
If $\lambda\notin\sigma_b(\m)$, then
$\lambda\in\Phi(A)\cap\Phi(B)$,
$\asc(A-\lambda)=\dsc(A-\lambda)<\infty$ and
$\asc(B-\lambda)=\dsc(B-\lambda)<\infty$. Since
$\dsc(A-\lambda)\Longrightarrow \lambda\notin\A$ and
$\asc(B-\lambda)\Longrightarrow \lambda\notin\b$, $\lambda\notin
\A\cap\b$. Hence, since $\sigma_b(\M)\subseteq\sigma_b(\m)$,
$\lambda\notin\sigma_b(\M)\cup\{\A\cap\b\}$. Conversely, if
$\lambda\notin\sigma_b(\M)\cup\{\A\cap\b\}$, then
$\lambda\in\Phi_+(A)\cap\Phi_-(B)$, $\asc(A-\lambda<\infty$
($\Longrightarrow \ind(A-\lambda)\leq 0$),
$\dsc(B-\lambda)<\infty$ ($\Longrightarrow \ind(B-\lambda)\geq
0$) and $\ind(A-\lambda)+\ind(B-\lambda)=0$. If $A^*$ has SVEP,
then $\ind(A-\lambda)\geq 0$; hence $\ind(A-\lambda)=0$, which
implies that $\ind(B-\lambda)=0$. But then both $A-\lambda$ and
$B-\lambda$ have finite (hence, equal) ascent and descent. Thus
$\lambda\notin\sigma_b(\m)$. Arguing similarly in the case in
which $\lambda\notin\b$ (this time using the fact that
$\lambda\in\Phi_-(B)$, $\dsc(B-\lambda)<\infty$ and
$\lambda\notin\b$ imply $\ind(B-\lambda)=0$), it is seen (once
again) that $\lambda\notin\sigma_b(\m)$.\end{proof} If we let
$$\sigma_{sb}(T)=\{\lambda\in\C:\hspace{2mm}\mbox{either}
\hspace{2mm}\lambda\notin\Phi_-(T)\hspace{2mm}\mbox{or}\hspace{2mm}
\dsc(T-\lambda)=\infty\},$$ then
$\sigma_{sb}(T)=\sigma_{ab}(T^*)$,
$\sigma_b(T)=\sigma_{ab}(T)\cup\sigma_{sb}(T)$ and
$\sigma_b(T)=\sigma_{ab}(T)\cup\Xi(T^*)=\sigma_{sb}(T)\cup\Xi(T)$
\cite[p 141]{A}. Evidently,
$\sigma_{ab}(\m)\cup\{\A\cup\B\}=\sigma_b(\m)$. Let $\Xi_+(T)$ and
$\Xi^*_+(T)$ denote the sets of $\lambda$ such that
$$\lambda\notin\Xi_+(T)\Longrightarrow \lambda\in\Phi_+(T)
\hspace{2mm}\mbox{and}\hspace{2mm} T\hspace{2mm}\mbox{has SVEP
at}\hspace{2mm}\lambda$$ and
$$\lambda\notin\Xi^*_+(T)\Longrightarrow \lambda\in\Phi_+(T)\hspace{2mm}\mbox{and}\hspace{2mm}
T^*\hspace{2mm}\mbox{has SVEP
at}\hspace{2mm}\lambda.$$\begin{prop}\label{prop2}
$\sigma_{ab}(\M)\subseteq\sigma_{ab}(\m)\subseteq\sigma_{ab}(\M)\cup\{\Xi^*_+(A)\cup\Xi_+(B)\}$.\end{prop}
\begin{proof} The inclusion $\sigma_{ab}(\M)\subseteq\sigma_{ab}(\m)$ being
evident, we prove
$\sigma_{ab}(\m)\subseteq\sigma_{ab}(\M)\cup\Xi^*_+(A)$. Let
$\lambda\notin\sigma_{ab}(\M)\cup\Xi^*_+(A)$. Then
$\lambda\in\Phi_+(A)$, $\asc(A-\lambda)<\infty$, and $A^*$ has
SVEP at $\lambda$. Hence $\ind(A-\lambda)=0$ and
$\lambda\in\Phi(A)$ ($\Longrightarrow
\lambda\notin\sigma_{ab}(A)$), which implies that
$\lambda\in\Phi_+(B)$ and $\ind(B-\lambda)\leq 0$(since $\asc(\M
-\lambda)<\infty\Longrightarrow
\ind(A-\lambda)+\ind(B-\lambda)\leq 0$). But then the hypothesis
$\lambda\notin\Xi_+(B)$ implies that
$\lambda\notin\sigma_{ab}(B)\Longrightarrow
\lambda\notin\sigma_{ab}(\m)$.\end{proof} The following corollary
is immediate from Proposition \ref{prop2}. \begin{cor}\label{cor1}
If $\Xi_+^*(A)\cup\Xi_+(B)=\emptyset$, then
$\sigma_{ab}(\m)=\sigma_{ab}(\M)$.\end{cor} It is easy to see
(from the definitions of $\sigma_b(T)$ and $\sigma_w(T)$) that
$$\sigma_b(T)=\sigma_w(T)\cup\Xi(T)=\sigma_w(T)\cup\Xi(T^*).$$
Hence
$$\sigma_b(\m)=\{\sigma_w(A)\cup\sigma_w(B)\}\cup\{\Xi(P)\cup\Xi(Q)\},$$
where $P=A$ or $A^*$ and $Q=B$ or $B^*$.\begin{prop}\label{prop3}
$\sigma_w(A)\cup\sigma_w(B)\subseteq\sigma_w(\M)\cup\{\Xi(P)\cup\Xi(Q)\}$,
where $P=A$ and $Q=B$ or $P=A^*$ and
$Q=B^*$.\end{prop}\begin{proof} The proof in both the cases is
similar: we consider $P=A$ and $Q=B$. If
$\lambda\notin\sigma_w(\M)$, then
$\lambda\in\Phi_+(A)\cap\Phi_-(B)$ and
$\ind(A-\lambda)+\ind(B-\lambda)=0$. Thus, since
$\lambda\notin\a\cup\b$ implies $\ind(A-\lambda)\leq 0$ and
$\ind(B-\lambda)\leq 0$, $\ind(A-\lambda)=\ind(B-\lambda)=0$ and
$\lambda\in\Phi(A)\cap\Phi(B)$. Hence
$\lambda\notin\sigma_w(A)\cup\sigma_w(B)$.\end{proof} Proposition
\ref{prop3} implies that if $\Xi(P)\cup\Xi(Q)=\emptyset$, $P$ and
$Q$ as above, then
$\sigma_w(\m)=\sigma_w(\M)=\sigma_w(A)\cup\sigma_w(B)$. More is
true. Since $\lambda\notin\sigma_w(\M)\cup\{\Xi(P)\cup\Xi(Q)\}$
implies $\lambda\in\Phi(A)\cap\Phi(B)$ and
$\ind(A-\lambda)=\ind(B-\lambda)=0$,
$\asc(A-\lambda)=\dsc(A-\lambda)<\infty$ and
$\asc(B-\lambda)=\dsc(B-\lambda)<\infty$. Hence
$\sigma_b(\M)\subseteq\sigma_w(\M)\cup\{\Xi(P)\cup\Xi(Q)\}$.
\begin{cor}\label{cor2} If $\Xi(P)\cup\Xi(Q)=\emptyset$, $P$ and $Q$ as
in Proposition \ref{prop3}, then
$\sigma_b(\m)=\sigma_w(\m)=\sigma_b(\M)=\sigma_w(\M)=\sigma_w(A)\cup\sigma_w(B)$.\end{cor}
\begin{proof}
$\sigma_w(\M)\subseteq\sigma_b(\M)\subseteq\sigma_b(\m)$.\end{proof}
The following theorem gives a necessary and sufficient condition
for $\sigma_w(\m)=\sigma_b(\m)$ and
$\sigma_{aw}(\m)=\sigma_{ab}(\m)$. \begin{thm}\label{thm1} (i)
$\sigma_w(\m)=\sigma_b(\m)$ if and only if $A$ and $B$ have SVEP
on $\{\lambda:\lambda\in\Phi(A)\cap\Phi(B),
\ind(A-\lambda)+\ind(B-\lambda)=0\}$.

(ii) $\sigma_{aw}(\m)=\sigma_{ab}(\m)$ if and only if $A$ and $B$
have SVEP on $\{\lambda:\lambda\in\Phi_+(A)\cap\Phi_+(B),
\ind(A-\lambda)+\ind(B-\lambda)\leq 0\}$.\end{thm}\begin{proof}
(i) If $\sigma_w(\m)=\sigma_b(\m)$, then
$\sigma_w(\m)=\sigma_w(A)\cup\sigma_w(B)=\sigma_b(A)\cup\sigma_b(B)$.
Hence $\lambda\in\Phi(A)\cap\Phi(B)$ with
$\ind(A-\lambda)=\ind(B-\lambda)=0$ if and only if
$\lambda\in\Phi(A)\cap\Phi(B)$,
$\asc(A-\lambda)=\dsc(A-\lambda)<\infty$ and
$\asc(B-\lambda)=\dsc(B-\lambda)<\infty$. Evidently, $A$ and $B$
have SVEP at points $\{\lambda:\lambda\in\Phi(A)\cap\Phi(B),
\ind(A-\lambda)+\ind(B-\lambda)=0\}$. Conversely, if
$\lambda\notin\sigma_w(\m)$, then $\lambda\in\Phi(A)\cap\Phi(B)$
and $\ind(A-\lambda)+\ind(B-\lambda)=0$. Since $A$ and $B$ have
SVEP at $\lambda$, $\ind(A-\lambda)$ and $\ind(B-\lambda)$ are
both $\leq 0$. Hence $\ind(A-\lambda)=\ind(B-\lambda)=0$, which
(because of SVEP) implies that
$\asc(A-\lambda)=\dsc(A-\lambda)<\infty$ and
$\asc(B-\lambda)=\dsc(B-\lambda)<\infty$. Thus
$\lambda\notin\sigma_b(A)\cup\sigma_b(B)\Longrightarrow
\sigma_b(\m)\subseteq\sigma_w(\m)$. Since
$\sigma_w(\m)\subseteq\sigma_b(\m)$ always,
$\sigma_w(\m)=\sigma_b(\m)$.

(ii) Since
$\sigma_{aw}(\m)\subseteq\sigma_{aw}(A)\cup\sigma_{aw}(B)
\subseteq\sigma_{ab}(A)\cup\sigma_{ab}(B)=\sigma_{ab}(\m)$,
$\sigma_{aw}(\m)=\sigma_{ab}(\m)$ implies that
$\sigma_{aw}(\m)=\sigma_{ab}(A)\cup\sigma_{ab}(B)$.
Equivalently,\begin{eqnarray*} &  &
\{\lambda:\lambda\in\Phi_+(A)\cap\Phi_+(B),
\ind(A-\lambda)+\ind(B-\lambda)\leq 0\}\\ & = &
\{\lambda:\lambda\in\Phi_+(A)\cap\Phi_+(B),
\asc(A-\lambda)<\infty, \asc(B-\lambda)<\infty\}.\end{eqnarray*}
Hence $A$ and $B$ have SVEP on
$\{\lambda:\lambda\in\Phi_+(A)\cap\Phi_+(B),
\ind(A-\lambda)+\ind(B-\lambda)\leq 0\}$. Conversely, if
$\lambda\notin\sigma_{aw}(\m)$, then
$\lambda\in\Phi_+(A)\cap\Phi_+(B)$ and
$\ind(A-\lambda)+\ind(B-\lambda)\leq 0$. Since $A$ and $B$ have
SVEP at $\lambda$, both $\asc(A-\lambda)$ and $\asc(B-\lambda)$
are finite. Hence
$\lambda\notin\sigma_{ab}(A)\cup\sigma_{ab}(B)$, which implies
that $\lambda\notin\sigma_{ab}(\m)\Longrightarrow
\sigma_{ab}(\m)\subseteq\sigma_{aw}(\m)$. Since
$\sigma_{aw}(\m)\subseteq\sigma_{ab}(\m)$ always, the proof is
complete.\end{proof}

\section{ Browder, Weyl Theorems}
Translating Theorem \ref{thm1} to the terminology of Browder's
theorem, $Bt$, and $a$-Browder's theorem, $a-Bt$, we see that
{\em a necessary and sufficient condition for $\m$ to satisfy $Bt$
is that $A$ and $B$ have SVEP at points
$\lambda\notin\sigma_w(\m)$, and that a necessary and sufficient
condition for $\m$ to satisfy $a-Bt$ is that $A$ and $B$ have
SVEP at points $\lambda\notin\sigma_{aw}(\m)$}. The following
theorem relates $Bt$ (resp., $a-Bt$) for $\m$ to $Bt$ (resp.,
$a-Bt$) for $\M$. Let $\sigma_{SF_+}(T)$ (resp.,
$\sigma_{SF_-}(T)$)  denote the {\em upper semi--Fredholm
spectrum} (resp., the {\em lower semi--Fredholm spectrum}) of
$T$.\begin{thm}\label{thm2} (a). If either (i) $A$ has SVEP at
points $\lambda\in\sigma_w(\m)\setminus\sigma_{SF_+}(A)$ and $B$
has SVEP at points $\mu\in\sigma_w(\m)\setminus\sigma_{SF_-}(B)$,
or (ii) both $A$ and $A^*$ have SVEP at points
$\lambda\in\sigma_w(\m)\setminus\sigma_{SF_+}(A)$, or (iii) $A^*$
has SVEP at points
$\lambda\in\sigma_w(\m)\setminus\sigma_{SF_+}(A)$ and $B^*$ has
SVEP at points $\mu\in\sigma_w(\m)\setminus\sigma_{SF_-}(B)$, then
$\m$ satisfies $Bt$ implies $\M$ satisfies $Bt$.

(b) If either (i) $A$ has SVEP at points
$\lambda\in\sigma_{aw}(\m)\setminus \sigma_{SF_+}(A)$ and $A^*$
has SVEP at points $\mu\in\sigma_w(\m)\setminus\sigma_{SF_+}(A)$,
or (ii) $A^*$ has SVEP at points
$\lambda\in\sigma_w(\m)\setminus\sigma_{SF_+}(A)$ and $B^*$ has
SVEP at points $\mu\in\sigma_w(\m)\setminus\sigma_{SF_+}(B)$, then
$\m$ satisfies $a-Bt$ implies $\M$ satisfies
$a-Bt$.\end{thm}\begin{proof} (a) Recall that $\m$ satisfies $Bt$
if and only if $M_0^*$ satisfies $Bt$, and that
$\sigma_w(M_0^*)=\sigma_w(\m)$. Hence $A$, $B$, $A^*$ and $B^*$
have SVEP at points $\lambda\notin\sigma_w(\m)$. In view of this,
hypotheses (i), (ii) and (iii) imply respectively that either
(i)' $A$ has SVEP at $\lambda\in\Phi_+(A)$ and $B$ has SVEP at
$\mu\in\Phi_-(B)$, or (ii)' $A$ and $A^*$ have SVEP at
$\lambda\in\Phi_+(A)$, or (iii)' $A^*$ has SVEP at
$\lambda\in\Phi_+(A)$ and $B^*$ has SVEP at $\mu\in\Phi_-(B)$.
Evidently, $\sigma_w(\M)\subseteq\sigma_b(\M)$: we prove that
$\sigma_b(\M)\subseteq\sigma_w(\M)$. For this, let
$\lambda\notin\sigma_w(\M)$. Then
$\lambda\in\Phi_+(A)\cap\Phi_-(B)$ and
$\ind(A-\lambda)+\ind(B-\lambda)=0$. Since both $\ind(A-\lambda)$
and $\ind(B-\lambda)$ are $\leq 0$ if (i)' holds,
$\ind(A-\lambda)=0$ if (ii)' holds, and both $\ind(A-\lambda)$ and
$\ind(B-\lambda)$ are $\geq 0$ if (iii)' holds, we conclude that
$\ind(A-\lambda)=\ind(B-\lambda)=0$,
$\lambda\in\Phi(A)\cap\Phi(B)$. Furthermore, since $A$ and $B$
have SVEP at $\lambda$, $\asc(A-\lambda)=\dsc(A-\lambda)<\infty$
and $\asc(B-\lambda)=\dsc(B-\lambda)<\infty$. Hence
$\lambda\notin\sigma_b(A)\cup\sigma_b(B)\Longrightarrow
\lambda\notin\sigma_b(\M)$.

\

(b) Since $a-Bt$ implies $Bt$, $M_0^*$ satisfies $Bt$, which
implies that $A^*$ has SVEP at $\lambda\notin\sigma_w(\m)$. Hence
the hypothesis $A^*$ has SVEP at
$\lambda\in\sigma_w(\m)\setminus\sigma_{SF_+}(A)$ implies that
$A^*$ has SVEP at $\lambda\in\Phi_+(A)$. Since the hypothesis $\m$
satisfies $a-Bt$ implies that $A$ and $B$ have SVEP on
$\{\lambda:\lambda\in\Phi_+(A)\cap\Phi_+(B),
\ind(A-\lambda)+\ind(B-\lambda)\leq 0\}$, it follows from the
SVEP hypotheses of the statement that either (i)' both $A$ and
$A^*$ have SVEP at $\lambda\in\Phi_+(A)$ or (ii)' $A^*$ has SVEP
at $\lambda\in\Phi_+(A)$ and $B^*$ has SVEP at $\mu\in\Phi_+(B)$.
Evidently, $\sigma_{aw}(\M)\subseteq\sigma_{ab}(\M)$. For the
reverse inclusion, let $\lambda\notin\sigma_{aw}(\M)$. Then
$\lambda\in\Phi_+(A)$. If (i)' is satisfied, then $A$ and $A^*$
have SVEP at $\lambda \Longrightarrow \ind(A-\lambda)=0$, which
(because $\lambda\in\Phi_+(A)$) implies that $\lambda\in\Phi(A)$
and $\ind(A-\lambda)=0$; if (ii)' holds, then $A^*$ has SVEP at
$\lambda$ implies $\ind(A-\lambda)\geq 0\Longrightarrow
\lambda\in\Phi(A)$ with $\ind(A-\lambda)\geq 0$. In either case
it follows that $\lambda\in\Phi_+(B)$  and
$\ind(A-\lambda)+\ind(B-\lambda)\leq 0$. Hence
$\lambda\in\Phi(A)\cap\Phi(B)$, $\ind(A-\lambda)=0$ and
$\ind(B-\lambda)\leq 0$. Since both $A$ and $B$ have SVEP on
$\{\lambda:\lambda\in\Phi_+(A)\cap\Phi_+(B),
\ind(A-\lambda)+\ind(B-\lambda)\leq 0\}$, $\asc(A-\lambda)<\infty$
and $\asc(B-\lambda)<\infty$. Thus
$\lambda\notin\sigma_{ab}(\m)\Longrightarrow
\lambda\notin\sigma_{ab}(\M)$.\end{proof}\begin{rem}\label{rem1}
 We note, for future reference, that if $\m$ satisfies $Bt$
and either of the hypotheses (i) to (iii) of Theorem
\ref{thm2}(a) is satisfied, then $\sigma_w(\m)=\sigma_w(\M)$.
Furthermore, $\sigma(\M)=\sigma(\m)$, as the following argument
shows. If $\m$ satisfies $Bt$, and one of the hypotheses (i),
(ii) and (iii) is satisfied, then either $A^*$ has SVEP at
$\lambda\in\Phi_+(A)$ or $B$ has SVEP at $\mu\in\Phi_-(B)$. Since
$\lambda\notin\sigma(\M)$ implies $A-\lambda$ is left invertible
and $B-\lambda$ is right invertible, $A^*$ has SVEP at
$\lambda\in\Phi_+(A)$ if and only if $A-\lambda$ is onto and $B$
has SVEP at $\lambda\in\Phi_-(B)$ if and only if $B-\lambda$ is
injective \cite[Corollary 2.4]{A}. In either case, both
$A-\lambda$ and $B-\lambda$ are invertible, which implies that
$\lambda\notin\sigma(\m)\Longrightarrow
\sigma(\m)\subseteq\sigma(\M)$. Since
$\sigma(\M)\subseteq\sigma(\m)$, the equality of the spectra
follows.\end{rem}
\begin{cor}\label{cor3} (a) \cite[Proposition 4.1]{DZ} If
$\{\a\cap\B\}\cup\A=\emptyset$, then $\m$ satisfies $Bt$ (resp.,
$a-Bt$) implies $\M$ satisfies $Bt$ (resp., $a-Bt$).

\

(b) \cite[Theorem 3.2]{Cao} If either
$\sigma_{aw}(A)=\sigma_{SF_+}(B)$ or
$\sigma_{SF_-}(A)\cap\sigma_{SF_+}(B)=\emptyset$, then $\m$
satisfies $Bt$ (resp., $a-Bt$) implies $\M$ satisfies $Bt$
(resp., $a-Bt$).\end{cor}\begin{proof}  Observe that if
$\{\a\cap\B\}\cup\A=\emptyset$, then either $A$ and $A^*$ have
SVEP or $A^*$ and $B^*$ have SVEP. Hence the proof for (a)
follows from Theorem \ref{thm2}(b). Assume now that
$\sigma_{aw}(A)=\sigma_{SF_+}(B)$. If
$\lambda\notin\sigma_{aw}(\M)$, then $\lambda\in\Phi_+(A)$ and
either $\alpha(B-\lambda)<\infty$ and
$\ind(A-\lambda)+\ind(B-\lambda)\leq 0$, or
$\beta(A-\lambda)=\alpha(B-\lambda)=\infty$ and $(B-\lambda){\X}$
is closed, or $\beta(A-\lambda)=\infty$ and $(B-\lambda){\X}$ is
not closed. Observe that if
$\alpha(B-\lambda)=\beta(A-\lambda)=\infty$ or
$\beta(A-\lambda)=\infty$, then the hypothesis
$\lambda\in\Phi_+(A)$ with $\ind(A-\lambda)\leq 0
\Longleftrightarrow\lambda\in\Phi_+(B)\Longrightarrow
\alpha(B-\lambda)<\infty$ -- a contradiction. Hence
$\lambda\in\Phi_+(A)\cap\Phi_+(B)$ and
$\ind(A-\lambda)+\ind(B-\lambda)\leq 0$. Again, if
$\sigma_{SF_-}(A)\cap\sigma_{SF_+}(B)=\emptyset$, then
$\Phi_-(A)\cup\Phi_+(B)=\C$. If $\lambda\notin\sigma_{aw}(\M)$,
then $\lambda\in\Phi_+(A)\Longrightarrow \lambda\in\Phi(A)$,
which (see above) implies that $\lambda\in\Phi_+(A)\cap\Phi_+(B)$
and $\ind(A-\lambda)+\ind(B-\lambda)\leq 0$. Hence if either of
the hypotheses $\sigma_{aw}(A)=\sigma_{SF_+}(B)$ and
$\sigma_{SF_-}(A)\cap\sigma_{SF_+}(B)=\emptyset$ holds, then
$\sigma_{aw}(\M)=\sigma_{aw}(\m)$. A similar argument, this time
using the fact that $\lambda\notin\sigma_w(\M)\Longrightarrow
\lambda\notin\sigma_{aw}(\M)=\sigma_{aw}(\m)$, shows that
$\sigma_w(\M)=\sigma_w(\m)$. (See \cite[Corollary 2.2 and Theorem
3.2]{Cao} for a slightly different argument.) Thus if $\m$
satisfies $Bt$ (resp., $a-Bt$), then $\M$ has SVEP at
$\lambda\notin\sigma_w(\M)$ (resp.,
$\lambda\notin\sigma_{aw}(\M)$), which implies that $\M$
satisfies $Bt$ (resp., $a-Bt$).\end{proof}\begin{rem}\label{rem2}
 If $\A\cup\B=\emptyset$, then $M_C^*$ has SVEP: this follows
from a straightforward application of the definition of SVEP
(applied to $(M_C^*-\lambda I^*)(f_1(\lambda)\oplus
f_2(\lambda))=0$). Hence $\sigma(\m)=\sigma(\M)=\sigma_a(\M)$,
$\sigma_{aw}(\M)=\sigma_w(\M)=\sigma_w(\m)$ and
$p_0(\M)=p_0^a(\M)$. Evidently, both $\m$ and $\M$ satisfy
$a-Bt$.\end{rem} We call an operator $T\in B(\Y)$ {\em polaroid}
\cite {DHJ} (resp., {\em isoloid}) at $\lambda\in\iso\sigma(T)$ if
$\asc(T-\lambda)=\dsc(T-\lambda)<\infty$ (resp., $\lambda$ is an
eigenvalue of $T$). Trivially, $T$ polaroid at $\lambda$ implies
$T$ isoloid at $\lambda$. Since
$$\pi_0(\m)=\{\pi_0(A)\cap\rho(B)\}\cup\{\rho(A)\cap\pi_0(B)\}\cup\{\pi_0(A)\cap\pi_0(B\},$$
if $\m$ is polaroid at $\lambda\in\pi_0(\m)$, then either $A$ or
$B$ is polaroid at $\lambda$; in particular, $A$  and $B$ are
polaroid at $\lambda\in\pi_0(A)\cap\pi_0(B)$. Conversely, if $A$
is polaroid at $\lambda\in\pi_0(A)$ and $B$ is polaroid  at
$\mu\in\pi_0(B)$, then $\m$ is polaroid at $\nu\in\pi_0(\m)$. We
say that $T$ is $a$-polaroid if $T$ is polaroid at
$\lambda\in\iso\sigma_a(T)$.
\begin{prop}\label{prop4} (i) $\m$ satisfies $Wt$ if and only if
$\m$ has SVEP at $\lambda\notin\sigma_w(\m)$ and $\m$ is polaroid
at $\mu\in\pi_0(\m)$.

(ii) $\m$ satisfies $a-Wt$ if and only if $\m$ has SVEP at
$\lambda\notin\sigma_{aw}(\m)$ and $\m$ is polaroid at
$\mu\in\pi_0^a(\m)$.\end{prop}\begin{proof} (i) is proved in
\cite[Theorem 2.2(i) and (ii)]{D}. To prove (ii) we start by
observing that if $\m$ has SVEP at
$\lambda\notin\sigma_{aw}(\m)$, then ($\m$ satisfies
$a-Bt\Longrightarrow$)
$\sigma_a(\m)\setminus\sigma_{aw}(\m)=p_0^a(\m)\subseteq\pi_0^a(\m)$,
which if points in $\pi_0^a(\m)$ are poles implies that
$\pi_0^a(\m)\subseteq p_0^a(\m)$. Conversely, $\m$ satisfies
$a-Wt$ implies $\m$ satisfies $a-Bt$, which in turn implies that
$\m$ has SVEP at $\lambda\notin\sigma_{aw}(\m)$. Again, since $\m$
(satisfies $a-Bt$ and) $a-Wt$, $\pi_0^a(\m)=p_0^a(\m)$.\end{proof}
A similar argument proves the following:
\begin{prop}\label{prop5} (i) $\M$ satisfies $Wt$ if and only if
$\M$ has SVEP at $\lambda\notin\sigma_w(\M)$ and $\M$ is polaroid
at $\mu\in\pi_0(\M)$.

(ii) $\M$ satisfies $a-Wt$ if and only if $\M$ has SVEP at
$\lambda\notin\sigma_{aw}(\M)$ and $\M$ is polaroid at
$\mu\in\pi_0^a(\M)$.\end{prop} The following theorem gives a
necessary and sufficient condition for $\M$ to satisfy $Wt$ in
the case in which either of the hypotheses (i), (ii) and (iii) of
Theorem \ref{thm2} is satisfied.
\begin{thm}\label{thm3} If either of the SVEP
hypotheses (i), (ii) and (iii) of Theorem \ref{thm2}(a) is
satisfied, then $\M$ satisfies $Wt$ for every $C\in B(\X)$ if and
only if $\m$ satisfies $Wt$ and $A$ is polaroid at
$\lambda\in\pi_0(\M)$.\end{thm}\begin{proof} {\em Sufficiency.} If
$\m$ satisfies $Wt$ (hence, $Bt$) and either of the hypotheses
(i), (ii) and (iii) of Theorem \ref{thm2} is satisfied, then $A$
satisfies $Bt$, $\sigma(\M)=\sigma(\m)$,
$\sigma_w(\M)=\sigma_w(\m)$ and $\M$ satisfies $Bt$ (see Remark
\ref{rem1} and Theorem \ref{thm2}(a)). Hence
$$\sigma(\M)\setminus\sigma_w(\M)=\sigma(\m)\setminus\sigma_w(\m)
=p_0(\m)=\pi_0(\m)\subseteq\pi_0(\M),$$ where the final inclusion
follows from the fact that
$\sigma(\M)\setminus\sigma_w(\M)=p_0(\M)\subseteq\pi_0(\M)$.
Hence to prove sufficiency, we have to prove the reverse
inclusion. Let $\lambda\in\pi_0(\M)$. Then
$\lambda\in\iso\sigma(\m)$. Start by observing that
$(\M-\lambda)^{-1}(0)\neq\emptyset\Longrightarrow
(\m-\lambda)^{-1}(0)\neq\emptyset$; also,
$\dim(\M-\lambda)^{-1}(0)<\infty\Longrightarrow
\dim(A-\lambda)^{-1}(0)<\infty$. We claim that
$\dim(B-\lambda)^{-1}(0)<\infty$. For suppose to the contrary
that $\dim(B-\lambda)^{-1}(0)$ is infinite. Since
$$(\M-\lambda)(x\oplus y)=\{(A-\lambda)x+Cy\}\oplus
(B-\lambda)y,$$ either $\dim(C(B-\lambda)^{-1}(0))<\infty$ or
$\dim(C(B-\lambda)^{-1}(0))=\infty$. If
$\dim(C(B-\lambda)^{-1}(0))<\infty$, then $(B-\lambda)^{-1}(0)$
contains an orthonormal sequence $\{y_j\}$ such that
$(\M-\lambda)(0\oplus y_j)=0$ for all $j=1,2,...$. But then
$\dim(\M-\lambda)^{-1}(0)=\infty$, a contradiction. Assume now
that $\dim(C(B-\lambda)^{-1}(0))=\infty$. Since $\lambda\in
\rho(A)\cup\iso\sigma(A)$, $A$ satisfies $Bt$, $A$ is polaroid
 at $\lambda\in\pi_0(\M)$ and $\alpha(A-\lambda)<\infty$,
$\beta(A-\lambda)<\infty$. Hence $\dim\{C(B-\lambda)^{-1}(0)\cap
(A-\lambda)\X\}=\infty$ implies the existence of a sequence
$\{x_j\}$ such that $(A-\lambda)x_j=Cy_j$ for all $j=1,2,...$ .
But then $(\M-\lambda)(x_j\oplus -y_j)=0$ for all $j=1,2,...$ .
Thus $\dim(\M-\lambda)^{-1}(0)=\infty$, again a contradiction.
Our claim having been proved, we conclude that
$\lambda\in\pi_0(\m)$. Thus $\pi_0(\M)\subseteq\pi_0(\m)$.

{\em Necessity.} Evidently, $\M$ satisfies $Wt$ for all $C$
implies $\m$ satisfies $Wt$. Hence
$p_0(\M)=\pi_0(\M)=p_0(\m)=\pi_0(\m)$, which implies that $\m$ is
polaroid at points $\lambda\in\pi_0(\M)$. Since
$\pi_0(\M)=p_0(\m)$, and since $\lambda\in p_0(\m)$ implies
$\lambda\in p_0(A)\cup\rho(A)$, $A$ is polaroid at
$\lambda\in\pi_0(\M)$.\end{proof}  \begin{rem}\label{rem3} An
examination of the proof of the sufficiency part of the theorem
above shows that if either of the SVEP hypotheses (i), (ii) and
(iii) of Theorem \ref{thm2}(a) is satisfied and $\m$ satisfies
$Wt$, then either of the hypotheses that $A$ is polaroid or $A$
is isoloid and satisfies $Wt$ is sufficient for $\M$ to satisfy
$Wt$.\end{rem}

\begin{cor}\label{cor4} (a) \cite[Theorem 4.2]{DZ} If
$\{\a\cap\B)\cup\A=\emptyset$, $A$ is polaroid at
$\lambda\in\pi_0(\M)$ (or $A$ is isoloid and satisfies $Wt$) and
$\m$ satisfies $Wt$, then $\M$ satisfies $Wt$.

(b) \cite[Theorem 3.3]{Cao} If $\sigma_{aw}(A)=\sigma_{SF_+}(B)$
or $\sigma_{SF_-}(A)\cap\sigma_{SF_+}(B)=\emptyset$, $A$ is
polaroid at $\lambda\in\pi_0(\M)$ (or $A$ is isoloid and
satisfies $Wt$) and $\m$ satisfies $Wt$, then $\M$ satisfies
$Wt$.\end{cor}\begin{proof} (a) Theorem \ref{thm3}, and Remark
\ref{rem3}, apply.\\ (b) Recall, \ref{cor3}(b),  that $\M$
satisfies $Bt \Longleftrightarrow
\sigma(\M)\setminus\sigma_w(\M)=p_0(\M)$. Hence
$\sigma(\M)\setminus\sigma_w(\M)\subseteq\pi_0(\M)$. For the
reverse inequality, start by recalling from the proof of
Corollary \ref{cor3}(b) that $\sigma_w(\M)=\sigma_w(\m)$. If
$\lambda_0\in\pi_0(\M)$, then there exists an $\epsilon$--
neighbourhood ${\mathcal N}_\epsilon$ of $\lambda_0$ such that
$\M-\lambda$ is invertible (implies $A-\lambda$ is left
invertible and $B-\lambda$ is right invertible), hence Weyl, for
all $\lambda\in{\mathcal N}_\epsilon$ not equal to $\lambda_0$.
Thus $\m-\lambda$ is Weyl for all $\lambda\in{\mathcal
N}_\epsilon$ not equal to $\lambda_0$. Since $\m$ satisfies $Bt$,
$\m-\lambda$ is Browder for all $\lambda\in{\mathcal N}_\epsilon$
not equal to $\lambda_0$, which implies that both $A-\lambda$ and
$B-\lambda$ are invertible. Hence $\lambda\in\iso\sigma(\m)$. Now
argue as in the sufficiency part of the proof of Theorem
\ref{thm3}.\end{proof}

The following examples, \cite{Lee} and \cite{DZ}, show that $\M$
in theorem above may fail to satisfy $Wt$ if one assumes only
that $A$ is isoloid but not polaroid at $\lambda\in\pi_0(\M)$, or
only that $A$ is polaroid at
$\lambda\in\pi_0(A)$.\begin{ex}\label{ex1}  Let $A$, $B$ and
$C\in B(\ell^2)$ be the operators
\begin{eqnarray*} &  & A(x_1,x_2,x_3,...)=(0,x_1,0,\frac{1}{2}
x_2,0,\frac{1}{3} x_3,...),\\ &  &
B(x_1,x_2,x_3,...)=(0,x_2,0,x_4,0,...),\end{eqnarray*} and
\begin{eqnarray*}
C(x_1,x_2,x_3,...)=(0,0,x_2,0,x_3,...).\end{eqnarray*} Then
$A$,$A^*$, $B$ and $B^*$ have SVEP, $\sigma(A)=\sigma_w(A)=\{0\}$,
$\pi_0(A)=p_0(A)=\emptyset$, and $A$ satisfies Weyl's theorem.
Since $\sigma(\m)=\sigma_w(\m)=\{0,1\}$ and
$\pi_0(\m)=p_0(\m)=\emptyset$, $\m$ satisfies $Wt$. However,
since $\sigma(\M)=\sigma_w(\M)=\{0,1\}$ and $\pi_0(\M)=\{0\}$,
$\M$ does not satisfy $Wt$. Observe that $A$ is not polaroid on
$\pi_0(\M)$.

Again, let $A$, $B$ and $C\in B(\ell^2)$ be the operators

\begin{eqnarray*} &  &
A(x_1,x_2,x_3,...)=(0,0,0,\frac{1}{2} x_2,0,\frac{1}{3}
x_3,...),\\ &  &
B(x_1,x_2,x_3,...)=(0,x_2,0,x_4,0,...),\end{eqnarray*} and
\begin{eqnarray*}
C(x_1,x_2,x_3,...)=(x_1,0,x_2,0,x_3,...).\end{eqnarray*}Then $A$,
$B$ (and $C$) have SVEP, $\sigma(A)=\sigma_w(A)=\pi_0(A)=\{0\}$,
and $\sigma(B)=\sigma_w(B)=\{0,1\}, \pi_0(B)=p_0(B)=\emptyset$.
Since
$$\sigma(\m)=\sigma_w(\m)=\{0\}
\hspace{2mm}\mbox{and}\hspace{2mm} \pi_0(\m)=p_0(\m)=\emptyset,$$
$\m$ satisfies $Wt$. However, since
$$\sigma(\M)=\sigma_w(\M)=\{0,1\}\hspace{2mm}\mbox{and}\hspace{2mm}\pi_0(\M)=\{0\},$$
$\M$ does not satisfy $Wt$. Observe that $0\notin p_0(A)$;
 $A$  satisfies $Bt$,
but does not satisfy $Wt$.\end{ex}  More can be said in the case
in which $\A\cup\B=\emptyset$. Recall from Remark \ref{rem2} that
if $\A\cup\B=\emptyset$, then $M_C^*$ has SVEP and $\M$ satisfies
$a-Bt$.\begin{thm}\label{thm4} If $\A\cup\B=\emptyset$, $A$ is
polaroid at $\lambda\in\pi_0^a(\M)$ (or, $A$ is isoloid and
satisfies $Wt$) and $B$ is polaroid at $\mu\in\pi_0^a(B)$, then
$\M$ satisfies $a-Wt$.\end{thm}\begin{proof}  Since $A^*$ and
$B^*$ have SVEP, both $M_0^*$ and $M_C^*$ have SVEP. Hence  $M_C$
(also, $M_0$) satisfies $Bt$, which implies that
$\sigma(\M)\setminus\sigma_w(\M)=p_0(\M)\subseteq\pi_0(\M)$.
Apparently, $\sigma(\m)=\sigma(\M)=\sigma_a(\M)$,
$\sigma_w(\m)=\sigma_w(\M)=\sigma_{aw}(\M)$,
$\pi_0(\M)=\pi_0^a(\M)$ and $\iso\sigma(\M)=\iso\sigma(\m)$.
Following (part of) the argument of the proof of the sufficiency
part of Theorem \ref{thm3}, it follows that if
$\lambda\in\pi_0(\M)$, then $\lambda\in\pi_0(A)\cap\pi_0(B)$. By
assumption, both $A$ and $B$ are polaroid at $\lambda$. Hence
$\m$ is polaroid at $\lambda$, which implies that $\lambda\in
p_0(\m)$. Since $\m$ satisfies $Bt$,
$\lambda\notin\sigma_w(\m)=\sigma_w(\M)$, which in view of the
fact that $\M$ satisfies $Bt$ implies that $\lambda\in p_0(\M)$.
Hence $\sigma(\M)\setminus\sigma_w(\M)=\pi_0(\M)\Longrightarrow
\sigma_a(\M)\setminus\sigma_{aw}(\M)=\pi_0^a(\M)$, i.e., $\M$
satisfies $a-Wt$.\end{proof} Theorem \ref{thm4} holds for polaroid
operators $A$ and $B$: for the polaroid hypothesis implies that
$\m$ is polaroid, hence satisfies $Wt$, which by Theorem
\ref{thm3} implies that $\M$ satisfies $Wt$. If the operators $A$
and $B$ have SVEP, then $\m$ and $\M$ have SVEP,
$\sigma(\m)=\sigma(\M)=\sigma(M_C^*)=\sigma_a(M_C^*)$,
$\iso\sigma(M_0^*)=\iso\sigma(M_C^*)=\iso\sigma_a(M_C^*)$,
$\pi_0(M_C^*)=\pi_0^a(M_C^*)$ and
$\sigma_w(\m)=\sigma_w(\M)=\sigma_w(M_C^*)=\sigma_{aw}(M_C^*)$.
Evidently, $A^*$, $B^*$, $M_0^*$ and $M_C^*$ satisfy $Bt$; in
particular, $p_0(M_0^*)=p_0(M_C^*)\subseteq\pi_0(M_C^*)$.
\begin{cor}\label{cor5} If the polaroid operators $A$ and $B$ have
SVEP, then $\M$ satisfies $Wt$ and $M_C^*$ satisfies
$a-Wt$.\end{cor}\begin{proof} Apparently, $\M$ satisfies $Wt$.
Since the polaroid hypothesis on $A$ and $B$ implies that $A^*$
and $B^*$ are polaroid, an argument similar to that in the
theorem above applied to $M_C^*$ implies that if
$\lambda\in\pi_0(M_C^*)$, then
$\lambda\in\pi_0(B^*)\cap\pi_0(A^*)\Longrightarrow \lambda\in
p_0(B^*)\cap p_0(A^*)\Longrightarrow
\lambda\notin\sigma_w(M_0^*)=\sigma_w(M_C^*)\Longrightarrow
M_C^*$ satisfies $Wt$. Hence $M_C^*$ satisfies $a-Wt$.\end{proof}

\bigskip

\noindent 8 Redwood Grove\\
Ealing\\
London W5 4SZ\\
United Kingdom

\

\noindent email: bpduggal@yahoo.co.uk

\end{document}